\documentclass[12pt]{article} 
\usepackage{latexsym}  
\usepackage{amsfonts} 
\usepackage{amsmath} 
\usepackage{amssymb} 
\usepackage{amscd}   
\usepackage{multicol}    
\usepackage{enumerate} 
\usepackage{amsthm}             
\theoremstyle{plain}
    \newtheorem{theorem}                    {Theorem}       [section]
    \newtheorem{lemma}      [theorem]       {Lemma}
    
    \newtheorem{proposition}[theorem]       {Proposition}

\begin{document}
\newcommand{\cosk}{\operatorname{cosk}}
\newcommand{\rem}{{\noindent \bf Remark: }} 
\newcommand{\chr}{\operatorname{char}}
\newcommand{\Alb}{\operatorname{Alb}}
\newcommand{\Jac}{\operatorname{Jac}}
\newcommand{\Aut}{\operatorname{Aut}}
\newcommand{\Hom}{\operatorname{Hom}}
\newcommand{\Tor}{\operatorname{Tor}}
\newcommand{\End}{\operatorname{End}}
\newcommand{\Ext}{\operatorname{Ext}}
\newcommand{\Sym}{\operatorname{Sym}}
\newcommand{\Gal}{\operatorname{Gal}}
\newcommand{\Pic}{\operatorname{Pic}}
\newcommand{\ord}{\operatorname{ord}}
\newcommand{\red}{\operatorname{red}}
\newcommand{\alb}{\operatorname{alb}}
\newcommand{\tor}{{\operatorname{tor}}}
\newcommand{\Spec}{\operatorname{Spec}}
\newcommand{\trdeg}{\operatorname{trdeg}}
\newcommand{\holim}{\operatornamewithlimits{holim}}
\newcommand{\im}{\operatorname{im}}
\newcommand{\coim}{\operatorname{coim}}
\newcommand{\coker}{\operatorname{coker}}
\newcommand{\gr}{\operatorname{gr}}
\newcommand{\id}{\operatorname{id}}
\newcommand{\Br}{\operatorname{Br}}
\newcommand{\cd}{\operatorname{cd}}
\newcommand{\CH}{\operatorname{CH}}
\renewcommand{\lim}{\operatornamewithlimits{lim}}
\newcommand{\colim}{\operatornamewithlimits{colim}}
\newcommand{\rk}{\operatorname{rank}}
\newcommand{\codim}{\operatorname{codim}}
\newcommand{\NS}{\operatorname{NS}}
\newcommand{\N}{{\mathbb N}}
\newcommand{\Z}{{{\mathbb Z}}}
\newcommand{\zl}{{{(l)}}}
\newcommand{\Q}{{{\mathbb Q}}}
\newcommand{\R}{{{\mathbb R}}}
\newcommand{\F}{{{\mathbb F}}}
\newcommand{\G}{{{\mathbb G}}}
\newcommand{\OO}{{{\mathcal O}}}
\renewcommand{\L}{{\mathcal L}}
\newcommand{\A}{{\mathcal A}}
\newcommand{\Sm}{\text{\rm Sm}}
\newcommand{\Sch}{\text{\rm Sch}}
\newcommand{\et}{{\text{\rm\'et}}}
\newcommand{\Zar}{{\text{\rm Zar}}}
\newcommand{\Nis}{{\text{\rm Nis}}}
\newcommand{\tr}{\operatorname{tr}}
\newcommand{\Div}{\operatorname{Div}}
\renewcommand{\div}{\operatorname{div}}
\newcommand{\corank}{\operatorname{corank}}
\renewcommand{\O}{{\cal O}}
\newcommand{\p}{{\mathfrak p}}
\newcommand{\proofend}{\hfill $\square$\\ \medskip}
\newcommand{\pschemes}{{\operatorname{Sch}^\bullet_k}}

\title{Rojtman's theorem for normal schemes}
\author{Thomas Geisser\thanks{Supported in part by JSPS grant-in-aid (B) 23340004 }}

\date{}
\maketitle

\begin{abstract}
We show that Rojtman's theorem holds for normal schemes: For every
reduced normal scheme $X$ of finite type over an algebraically closed
field $k$, the torsion subgroup of the zero'th Suslin homology 
is isomorphic to the torsion subgroup of
the $k$-rational points of the albanese variety of $X$
(the universal object for morphisms to semi-abelian varieties).
\end{abstract}

\section{Introduction} 
By a classical theorem of Abel and Jacobi, there is an isomorphism 
between the Chow group of zero-cycles of degree $0$ and the rational points of the Jacobian 
variety  
$$ CH_0(C)^0 \stackrel{\sim}{\longrightarrow} \Jac_C(k)$$
for every smooth projective curve $C$ over an algebraically closed field.
Rosenlicht gave a generalization to smooth (open) curves, 
comparing a Chow group with modulus to the generalized Jacobian variety
\cite{ros}, an extension of an abelian variety by a torus.

If $X$ is a smooth and projective scheme of 
higher dimension, it is natural generalization is to study the albanese map 
$$\alb_X:  CH_0(X)^0 \to \Alb_X(k) .$$
It is surjective, but can have a large kernel if 
$X$ has dimension at least $2$ \cite{mumford}. However,
Rojtman \cite{rojtman} proved that $\alb_X$ induces an isomorphism of torsion subgroups away from
the characteristic. A cohomological proof of Rojtman's theorem has been given by 
Bloch \cite{bloch}, and  Milne \cite{milnerojtman}  proved the analogous 
statement for the $p$-part in characteristic $p$.

Rojtman's theorem has been generalized in several directions. 
If $X$ is projective, then, using an improved duality theorem, the method of Bloch and Milne carries over to generalize Rojtman's theorem to normal schemes 
\cite{ichduality}. 
If $X$ is an open subscheme of a smooth projective scheme, then 
Spiess-Szamuely \cite{ss} showed that away from the characteristic, 
there is an isomorphism 
$$ \alb_X :{}_\tor H_0^S(X,\Z)^0\to {}_\tor \Alb_X(k).$$ 
Here the left hand side is Suslin homology, and the right hand side is 
Serre's albanese variety \cite{serre}, the universal semi-abelian
variety (i.e. extension of an abelian variety by a torus) to which $X$
maps. This theorem was generalized by Barbieri-Viale and Kahn 
\cite[Cor. 14.5.3]{bvk} to normal schemes in characteristic $0$. 
The theorem of Spiess-Szamuely holds for the $p$-part in characteristic $p$ 
under resolution of singularities \cite{ichsuslin}. In this paper, we
give a different proof of the theorem of Barbieri-Viale and Kahn 
which also gives the result in characteristic $p$, including
the $p$-part:

\begin{theorem} \label{rojtman}
Let $X$ be a reduced normal scheme, 
separated and of finite type over an algebraically closed field $k$
of characteristic $p\geq 0$. Then the albanese map
induces an isomorphism 
$${}_\tor H_0^S(X,\Z) \stackrel{\sim}{\longrightarrow} {}_\tor \Alb_X(k)$$
up to $p$-torsion groups, and 
$H_1^S(X,\Z)\otimes\Q_l/\Z_l=0$ for $l\not=p$.
Under resolution of singularities for schemes of dimension at most 
$\dim X$, the restriction on the characteristic
is unnecessary.
\end{theorem}

We note that there are also generalizations of Rojtman's theorem
comparing the cohomological Chow group (as defined by Levine-Weibel)
with Lang's albanese variety (universal for rational maps to abelian
varieties)
see \cite{bvpw}\cite{bs}\cite{krsr}\cite{mal}.

The idea of our proof 
is to work with the {\it albanese scheme} $\A_X$ introduced by Ramachandran.
It is universal for maps from $X$ to {\it locally semi-abelian schemes},
 i.e. group schemes 
locally of finite type whose connected component is a semi-abelian variety
and whose group of components is a lattice. 
The connected component of $\A_X$ is isomorphic to the albanese variety 
for every choice of a base-point of $X$. Using a  suggestion of Ramachandran, we prove

\begin{theorem}\label{albaba}
Let $X$ be a reduced, semi-normal, connected variety over a perfect field, 
and $a:X_\bullet\to X$ be a $1$-truncated proper hypercover $X$ such 
that $X_0\to X$ is proper and generically etale.
Then the albanese scheme $\A_X$ of $X$ is the largest locally semi-abelian scheme
quotient of $\A_{X_0}/d\A_{X_1}$, where $d=(\delta_0)_*-(\delta_1)_*$ for
$\delta_i:X_1\to X_0$ the two face maps.
\end{theorem}

On the other hand, for each prime $l$ different from the characteristic 
of $k$,
Suslin homology tensored with $\Z_{(l)}$ can be calculated by a proper 
ldh-covering \cite{ichdescent}
(and by a hyperenvelope at the characteristic). 
This allows us to prove the main theorem by reducing to the theorem of 
Spiess-Szamuely. 

\smallskip

If $X$ is not normal, then Rojtman's theorem is wrong even for curves,
as one sees by taking an elliptic curve and identifying $0$ with a non-torsion point
\cite{ichduality}. However, we propose a statement in terms of a hypercover which could 
serve as a generalization of Rojtman's theorem, and prove it for curves.
The statement for curves has he following explicit version:

\begin{theorem}
Let $C$ be a reduced semi-normal curve over an algebraically closed
field with normalization $\tilde C$. 
Then the albanese map induces an isomorphism 
$$H_0^S(C,\Z) \cong \Z[\pi_0(C)] \oplus \A_{\tilde C}^0(k) /H_1(D_{C_\bullet},\Z)$$
\end{theorem}

Here $H_1(D_{C_\bullet},\Z)$ is a free abelian group dual to $H^1_\et(C,\Z)$,
and $\A_{\tilde C}^0$ the connected component of the albanese scheme. 
Dividing the rational points $\A_{\tilde C}^0(k)$
by a free subgroup makes the torsion subgroup larger, explaining the example
in \cite{ichduality}.  

\medskip
The results of this paper were reported on in \cite{ichkyoto}. 

\medskip
Notation: For abelian group $A$, we write $A_{(l)}$ for
$A\otimes_\Z \Z_{(l)}$, $A[l]$  for the subgroup of $l$-power torsion 
elements, and ${}_\tor A=\oplus_l A[l]$ for the subgroup of torsion 
elements. 

\medskip
Acknowledgements: We thank N.Ramachandran for discussions on his
work, and the referee for his careful reading and helpful suggestions.

\section{The Albanese scheme}
Throughout this paper, $k$ is a perfect field, and $\Sch/k$ the
category of schemes locally of finite type over $k$.
The following discussion is based on work of Ramachandran \cite{ramachandran}, 
the notation follows Kahn-Sujatha \cite{ks}. 
A {\it semi-abelian variety} is an extension of an abelian variety by a torus,
and a {\it locally semi-abelian scheme} is a commutative group scheme such
that the scheme of components $\pi_0(\A)$ is a lattice $D$, and the connected component $\A^0$ is a
semi-abelian variety. 
For a homomorphism between semi-abelian varieties  $d:\A_1^0\to \A_0^0$,
the image of $d$ is a closed subvariety, the quotient $\A_1^0/d \A_0^0$
is a again a semi-abelian variety, and, if $k$ is algebraically closed,  
$\A_1^0(k)/d \A_0^0(k)$ is isomorphic to $\big(\A_1^0/d \A_0^0\big) (k)$.

\begin{lemma}\label{firstlemma}
Let $d:\A_1\to \A_0$  be a homomorphism of locally semi-abelian schemes 
and $\A$ be the largest locally semi-abelian scheme quotient of the presheaf
quotient $\A_0/d\A_1$. Then 
$ \A^0$ is the largest semi-abelian variety quotient of 
$$\A^0_{0} / \big( d\A_{1}\cap \A^0_{0}\big) \cong 
\A^0_{0} / \langle d\A_{1}^0, \delta \ker (D_1\to D_0)\rangle .$$
Here $\delta $ is the map from the snake Lemma in the following diagram
\begin{equation}\label{albmorph}
\begin{CD}
0@>>> \A_{1}^0 @>>> \A_{1} @>>> D_{1}@>>> 0\\
@. @VVV @VdVV @VVV\\
0@>>> \A_{0}^0 @>>> \A_{0} @>>> D_{0}@>>> 0\\
\end{CD}
\end{equation}
\end{lemma}

\proof
This follows from the definitions 
and the exact sequence of presheaves on $\Sch/k$:
$$\ker (D_1\to D_0)\stackrel{\delta}{\longrightarrow}
\A^0_{0} / d\A_{1}^0\to 
\A_{0}/d\A_{1}\to D_{0}/D_{1}\to 0.$$
\proofend

Consider the flat site $(Sch/k)_{fl}$ 
and $Z_X$ be the sheaf associated to 
the free abelian group on the presheaf $U\mapsto \Z[\Hom_k(U,X)]$
represented by $X$. The {\it albanese scheme} $u_X: X\to \A_X$
is the universal object for morphisms from $Z_X$ to sheaves represented by
locally semi-abelian schemes. 
For reduced schemes of finite type over $k$, the albanese scheme 
exists \cite[Thm.1.11]{ramachandran}, and the assignment $X\to \A_X$ 
is a covariant functor. 
By definition, there is an exact sequence
$$ 0\to \A_X^0\to \A_X \to D_X\to 0,$$
and $D_X$ is the flat sheaf associated
to the free presheaf $U\mapsto \Z\Hom(U,\pi_0(X))$, where  $\pi_0(X)$
is the largest quotient scheme of $X$ etale over $k$. 
For example, $D_X\cong \Z$ if $X$ is geometrically connected. 
For a connected scheme $X$, the scheme $\A_X^0$ is isomorphic to 
the usual Albanese variety $\Alb_X$ for every choice of a 
base-point $x_0$, because the map $X\to \A^0_X$, 
$x\mapsto u_X(x)-u_X(x_0)$ factors through $\Alb_X$ by
the universal property.

Recall that Suslin homology is the homology of the complex 
$C_*(X)$ which in degree $i$ consists of the free abelian group
generated by closed irreducible subschemes of $X\times \Delta^i$
which map finitely and surjectively onto $\Delta^i$. The boundary maps
are the alternating sums of pull-back maps to the faces. 
Let $H^S_0(X,\Z)^0$ be the kernel of the canonical degree map
(induced by covariant functoriality) 
$H^S_0(X,\Z) \to H^S_0(\pi_0(X),\Z)\cong D_X$.

\begin{lemma}
The albanese map induces a map from Suslin-homology to the 
albanese scheme such that the following diagram is commutative:
\begin{equation}
\begin{CD}
0@>>> H^S_0(X,\Z)^0  @>>> H^S_0(X,\Z)  @>>> D_{X}@>>> 0\\
@. @VVV @V\alb_X VV @|\\
0@>>> \A_{X}^0(k) @>>> \A_{X}(k) @>>> D_{X}@>>> 0\\
\end{CD}
\end{equation}
\end{lemma}

\proof
Extending the map $u_X$ to formal linear combinations, we obtain
a map from zero-cycles on $X$ to $\A_X$, and 
the Lemma follows exactly as in \cite[Lemma 3.1]{ss}.
\proofend

\subsection{Proper hypercovers}
Recall that a simplicial object and an $n$-truncated simplicial object
in a category $\mathcal C$ is a contravariant functor from the category of finite
totally ordered sets ${\mathbf \Delta}$ and the category ${\mathbf \Delta}^{\leq n}$ 
of finite totally ordered sets of order at most $n+1$
to $\mathcal C$, respectively. 
If $\mathcal C$ has finite limits, then the restriction functor
$i_n^*$ from simplicial sets to $n$-truncated simplicial sets has a 
right adjoint $(i_n)_*$, and we denote the composition $(i_n)_*i_n^*$
by $\cosk_n$.   
A proper hypercover $X_\bullet \to X$ is an augmented simplicial scheme 
$X_\bullet$ such that the adjunction maps 
$X_{i+1} \to (\cosk_i X_\bullet)_{i+1}$ 
are proper and surjective, and $n$-truncated proper hypercovers are defined similarly. 
For example, a $1$-truncated proper hypercover is a diagram
$$\begin{CD} 
X_1 \rightrightarrows  X_0 \to  X
\end{CD}$$
such that $a\delta_0=a\delta_1$, and such that $(\delta_0,\delta_1):X_1\to X_0\times_X X_0$ 
is proper and surjective,  together with a section $s:X_0\to X_1$
to $\delta_0$ and $\delta_1$.

{\it Proof of Theorem \ref{albaba}}.
Let $d={\delta_0}_*-{\delta_1}_*:\A_{X_1}\to \A_{X_0}$, and let 
$\A$ be the largest locally semi-abelian scheme quotient of 
$\A_{X_0}/d\A_{X_1}$.  Consider the following commutative diagram
$$\begin{CD}
X_1 @>\delta_0>\delta_1> X_0 @>a>> X\\
@Vu_1 VV @Vu_0VV @Vu_XVV \\
\A_{X_1} @>\delta_0>\delta_1> \A_{X_0} @>a>> \A_X,
\end{CD}$$
and denote the composition $X_0\stackrel{u_0}{\to}\A_{X_0}\to \A$ by $u'$.
Since $a\delta_0=a\delta_1$, the canonical map 
$\A_{X_0}\to \A_X$ factors through $\A_{X_0}/d\A_{X_1}$, 
hence through $\A$. It suffices to show that the induced map 
$$\A\to \A_X$$ 
is an isomorphism of locally semi-abelian schemes. 
Since $a$ is surjective, so
is the composition $\A_{X_0}\to \A\to \A_X$, and
it suffices to show that
the map $u_X:  X\to \A_X$ factors through $\A$.
Let $V\subseteq X$ be a dense open subset such that $V_0=a^{-1}V\to V$ is
etale and surjective, hence faithfully flat and thus a universal 
epimorphism, i.e. for any scheme $T$, the following sequence is
an equalizer:
$$ \Hom(V,T) \to \Hom(V_0,T)\rightrightarrows \Hom(V_0\times_VV_0,T).$$
Since $V_0\to V$ is etale and $V$ is reduced, 
$V_0\times_VV_0$ is reduced, hence after decreasing $V$ further, 
we can assume that the map 
$V_1=V\times_XX_1 \to V\times_X(X_0\times X_0)\cong V_0\times_VV_0$  induced by $(\delta_0,\delta_1)$ is (faithfully) flat.
Indeed being flat is an open condition, and since the target is reduced, 
surjectivity implies flatness at the generic point. 
Consequently $\Hom(V_0\times_VV_0,T)\to \Hom(V_1,T)$ is 
injective for any scheme $T$. 

Taking $T=\A$ and $u'|_{V_0}\in \Hom(V_0,\A)$, 
the two pull-back maps to $\Hom(V_1,\A)$, hence the pull-back maps to 
$\Hom(V_0\times_VV_0,\A)$ agree, so there is a unique map 
$v\in \Hom(V,\A)$ with $u'|_{V_0}= v\circ a|_{V_0}$.
It suffices to extend $v$ to $X$.
Consider the graphs $\Gamma_v\subset V\times \A$  
and $\Gamma_{u'}\subset X_0\times \A$ of $v$ and $u'$,
respectively, and let $C$ be the reduced
image of $\Gamma_{u'}$ in $X\times \A$. Consider the following diagram
$$\begin{CD}
\Gamma_v @>>> V\times \A @>>> V\\
@VinclVV @VinclVV @VinclVV \\
C @>>> X\times \A @>>> X\\
@A\alpha AA @Aa'AA @AaAA \\
\Gamma_{u'} @>>> X_0\times \A @>>> X_0.
\end{CD}$$
Since the upper and lower horizontal compositions are isomorphisms, 
$\alpha$ is proper (and surjective), and $C$ is closed in $X\times \A$.
It suffices to show that the middle composition $f:C\to X$ is an isomorphism.
Comparing to the upper row we see that $f$ is birational and proper, 
and by maximality of the 
semi-normalization it suffices to show that it is a bijection on $K$-rational
points for any field extension $k\subseteq K$. 
 
Let $x_1,x_2\in C(K)$ be two points with $f(x_1)=f(x_2)$, 
and $\tilde x_i$ be lifts to $\Gamma_{u'}\cong X_0$.
We claim that the image $u'(\tilde x_i)\in \A$ is independent of $i$.
Indeed, since $X_1\twoheadrightarrow X_0\times_X X_0$,
there is a point $t\in X_1$ such that $\delta_i(t)=\tilde x_i$, hence 
$u'(\tilde x_1)=u'(\tilde x_2)$. So 
$a'(\tilde x_1,u'(\tilde x_1))= a'(\tilde x_2,u'(\tilde x_2))$ and hence
$x_1=\alpha(\tilde x_1)=\alpha(\tilde x_2)=x_2$ because $C\to X\times \A$
is a closed embedding.
\proofend

\rem The statement of the theorem is wrong if $X$ is not semi-normal. 
For example, let $C$ be $\G_m$ with a cusp at the point $1$
and  $p:\G_m\to C$ be the normalization.
Then the simplicial albanese variety defined above is isomorphic to $\G_m$, 
but the (naive) albanese variety of $C$ is trivial. Indeed, being a quotient of $\G_m$, it
is either trivial or isomorphic to $\G_m$. If it was isomorphic to $\G_m$,
then from the commutativity of 
$$\begin{CD}
\G_m @>u_{\G_m} >\sim > \G_m\\
@VpVV @VVV \\
C@> u_C >> \G_m 
\end{CD}$$
we see that the composition
$u_Cp:\G_m\to \G_m$ is an isogeny of some degree $m$,
hence every closed point has $m$ inverse images. 
But as a dominant map to a regular curve, $u_C$ is flat, and over 
the origin, $p$ has degree $2$, a contradiction.

\section{Rojtman's theorem}
From now on we assume that our base field is algebraically closed.
Consider a $2$-truncated simplicial scheme $X_\bullet$ together with the
corresponding map of locally semi-abelian schemes 
\begin{equation}
\begin{CD}
0@>>> \A_{X_2}^0 @>>> \A_{X_2} @>>> D_{X_2}@>>> 0\\
@. @VVV @VdVV @VVV\\
0@>>> \A_{X_1}^0 @>>> \A_{X_1} @>>> D_{X_1}@>>> 0\\
@. @VVV @VdVV @VVV\\
0@>>> \A_{X_0}^0 @>>> \A_{X_0} @>>> D_{X_0}@>>> 0\\
\end{CD}
\end{equation}
with $d$ is the alternating sum of the maps induced by the
face maps $X_i\to X_{i-1}$. Taking $k$-rational points, we obtain
an analogous diagram of abelian groups. Let 
$H_i(\A_{X_\bullet}(k),\Q/\Z)= \Tor_i(\A_{X_\bullet}(k),\Q/\Z)$ 
be the homology of the double complex
$ \A_{X_\bullet}(k)\to \A_{X_\bullet}(k)_\Q$.
It is easy to see that $H_0(D_{X_\bullet},\Z)\cong \Z^{\pi_0(X)}$ is free,
so that the exact sequence of $k$-rational points
\begin{equation}\label{ssse}
H_1(D_{X_\bullet},\Z)  \stackrel{\delta}{\longrightarrow}
\A^0_{X_0}(k) / d\A_{X_1}^0(k)\to 
\A_{X_0}(k)/d\A_{X_1}(k)\to H_0(D_{X_\bullet},\Z)\to 0
\end{equation}
gives an isomorphism of abelian groups
$$ {}_\tor \big( \A^0_{X_0}(k)/ (d\A^0_{X_1}(k)+ \im \delta) \big)
\cong {}_\tor   \big( \A_{X_0}(k)/ d\A_{X_1} (k)\big).$$
We follow \cite{ramachandran} and \cite{bs} in defining the Albanese 
$1$-motive $M(X_\bullet)$ of $X_\bullet$ as
$$ H_1(D_{X_\bullet},\Z)/\tor \stackrel{\delta}{\to} 
\A^0_{X_0} / d\A_{X_1}^0.$$
Its homology with $\Q_l/\Z_l$-coefficients
$$H_1(M(X_\bullet),\Q_l/\Z_l):= 
H_1( M(X_\bullet)_\zl\to M(X_\bullet)_\Q)$$
sits in an exact sequence
\begin{equation}\label{onemotseq}
0\to (\A^0_{X_0}(k) /d \A_{X_1}^0(k))[l] \to H_1(M(X_\bullet),\Q_l/\Z_l)
\to H_1(D_{X_\bullet},\Z)\otimes\Q_l/\Z_l\to 0.
\end{equation}

\begin{theorem}\label{generalprop}
Let $X$ be separated and of finite type over the algebraically closed field $k$ of characteristic $p\geq 0$, 
and let $X_\bullet$ be a $2$-truncated proper hypercover of $X$ 
consisting of normal schemes. Then we have canonical isomorphisms
$$ H_1^S(X,\Q_l/\Z_l) \cong 
H_1(\A_{X_\bullet}(k),\Q_l/\Z_l)
\cong H_1(M(X_\bullet ),\Q_l/\Z_l)$$
if either $l\not=p$, or if resolution of singularities exists for
schemes of dimension at most $\dim X$.
\end{theorem}

\proof
Recall that for $l\not=p$, 
an $n$-truncated $l$-hyperenvelope (hyperenvelope) is a 
$n$-truncated proper hypercover $X_\bullet\to X$ satisfying the 
following condition:
for any $i\leq n$ and any point in the target of  
$X_{i+1} \to (\cosk_i X_\bullet)_{i+1}$,
there is a point mapping to it such that the extension of residue fields 
is finite of order prime to $l$ (trivial), respectively. We introduce 
the convention that for $l=p$, an $n$-truncated $l$-hyperenvelope 
means an $n$-truncated hyperenvelope. 

We first prove the theorem in case that $X_\bullet$ is a 
$2$-truncated proper $l$-hyperenvelope
of $X$ which is contained as an open subscheme in a 
$2$-truncated simplicial scheme $\bar X_\bullet$ consisting of smooth,
projective schemes. 

Let $X$ be a smooth scheme embedded in a smooth projective scheme
$\bar X$. Let $B$ be the image of $d_1$ in the Suslin complex 
$C_2(X)\stackrel{d_1}{\to} C_1(X)\stackrel{d_0}{\to} C_0(X)$, and consider the double 
complex $S(X)$ given by 
$$\begin{CD} 
(C_1(X)/B)_\zl/\tor  @>d_0>> C_0(X)_\zl \\
@VVV @VVV \\
(C_1(X)/B)_\Q @>d_0>> C_0(X)_\Q 
\end{CD} $$ 
Clearly $H_i(S(X))\cong H_i^S(X,\Q_l/\Z_l)$ for $i\leq 1$,
and since the left vertical map is injective, we obtain
$H_i(S(X))=0$ for $i\geq 2$.  

Let $T(X)$ be the complex $ \A_X(k)_\zl\to \A_X(k)_\Q$.
The albanese map induces a map of complexes $S(X) \to T(X)$, which we claim
to be a quasi-isomorphism. Indeed, in degree $0$ both groups are isomorphic
to  $D_X\otimes \Q_l/\Z_l $, and the map is compatible with 
this isomorphism by Lemma \ref{firstlemma}.
By our hypothesis on $X$, 
the theorem of Spiess-Szamuely \cite{ss} implies that 
$H_1^S(X,\Z)\otimes \Q_l/\Z_l=0$, hence in degree $1$ the map induces the 
isomorphism $H_0^S(X,\Z)[l]\cong \A_X(k)[l]$.

By the above discussion,  the rows of complexes coming from the hypercover,
$$\begin{CD}
S(X_2)@>>> S(X_1)@>>> S(X_0) \\
@VVV @VVV @VVV \\
T(X_2)@>>> T(X_1)@>>> T(X_0) .
\end{CD}$$
are quasi-isomorphic. 
By \cite{ichdescent}, the upper row calculates 
$H_1^S(X,\Q_l/\Z_l)$, whereas 
the lower row calculates $H_1(\A_{X_\bullet}(k),\Q_l/\Z_l)$
(here we use the fact that $X_\bullet$ is an $l$-hyperenvelope).

We claim that the canonical map 
$\alpha: M(X)\to \A_{X_0}(k)/d\A_{X_1}(k)$ induces the second 
isomorphism.
If we first take vertical homology in the double complex $T(X_\bullet)$, 
we obtain as the $E_1$-page 
$$ \begin{CD}
\A_{X_2}(k)[l] @>>> \A_{X_1}(k)[l] @>>> \A_{X_0}(k)[l]\\
@. @. @. \\
D_{X_2}\otimes \Q_l/\Z_l @>>> D_{X_1}\otimes \Q_l/\Z_l @>>> 
D_{X_0}\otimes \Q_l/\Z_l .
\end{CD}$$ 
Now 
$$ {}_\tor \A_{X_0}(k)/d {}_\tor \A_{X_1}(k) =
{}_\tor \A_{X_0}^0 (k)/d {}_\tor \A_{X_1}^0(k)\cong 
{}_\tor \big( \A_{X_0}^0(k)/d \A_{X_1}^0(k) \big) $$
on the one hand, and since $H_0(D_{X_\bullet},\Z)$ as well as the 
$D_{X_i}$ are torsion free, 
$$ H_1(D_{X_\bullet})\otimes \Q_l/\Z_l
\stackrel{\sim}{\to} H_1(D_{X_\bullet},\Q_l/\Z_l) 
=H_1(D_{X_\bullet}\otimes\Q_l/\Z_l) .$$
Hence  $\alpha$ induces a diagram of short exact sequences 
\begin{equation}\label{albaneseseq}
\begin{CD}
(\A^0_{X_0}(k) /d \A_{X_1}^0(k))[l] @>>> H_1(M(X_\bullet),\Q_l/\Z_l)
@>>> H_1(D_{X_\bullet},\Z)\otimes\Q_l/\Z_l\\
 @VVV @VVV @VVV\\
\A_{X_0}(k)[l] /d\A_{X_1}(k)[l] @>>> H_1(\A_{X_\bullet}(k),\Q_l/\Z_l)
@>>> H_1(D_{X_\bullet}\otimes\Q_l/\Z_l) .
\end{CD}\end{equation} 
in which the outer maps are isomorphisms, and hence so is the middle map.
This ends the proof of Theorem \ref{generalprop} in the case that 
$X_\bullet$ is a $2$-truncated $l$-hyperenvelope contained as an open
simplicial scheme in  a simplicial scheme 
$\bar X_\bullet$ consisting of smooth, projective schemes. 

\smallskip

By Gabber's refinement of de Jong's theorem on alterations \cite{gabber}
(or assuming resolution of singularities for schemes of dimenision at
most $\dim X$ if $l=p$), any scheme over
a perfect field admits an $l$-hyperenvelope $X_\bullet$ consisting of smooth schemes
which can be embedded into smooth projective schemes 
(see \cite[Thm.1.5]{ichdescent}). Hence 
for {\it any} reduced semi-normal scheme, the albanese map
induces a map of short exact coefficient sequences
\begin{equation}\label{ses3}
\begin{CD}
H_1^S(X)\otimes \Q_l/\Z_l @>>> H_1^S(X,\Q_l/\Z_l)@>>> H_0^S(X,\Z)[l] \\
@VVV @| @VVV \\
H_1(\A_{X_\bullet}(k))\otimes \Q_l/\Z_l @>>> 
H_1(\A_{X_\bullet}(k),\Q_l/\Z_l)@>>>(\A_{X_0}(k)/d\A_{X_1}(k))[l], 
\end{CD}
\end{equation}
and hence a surjection 
$$\alb_X: H_0^S(X,\Z)[l] \to  (\A_{X_0}(k)/d\A_{X_1}(k))[l]$$
if either $l\not=p$ or if resolution of singularities exists.
We next show that for {\it normal} $X$, the left hand terms of \eqref{ses3}
vanish, so that this surjection is an isomorphism.

Recall that 
$H^i_\et(X_\bullet,\mathcal F^\bullet)$ is the derived functor of the global section
functor $\mathcal F^\bullet \mapsto \ker \delta_0^*-\delta_1^* : 
\mathcal F^0(X_0)\to \mathcal F^1(X_1)$ from simplicial sheaves of 
abelian groups on $X_\bullet$ to abelian groups. 

\begin{proposition}
If $a: X_\bullet\to X$ is a $2$-truncated proper hypercover by normal schemes, then 
$$H^i_\et(X,\Z)\cong H^i_\et(X_\bullet,\Z)\cong
H^i(D_{X_\bullet},\Z)$$ for $i\leq 1$.
In particular,  $H_1(D_{X_\bullet},\Z)$ is finite for normal $X$.
\end{proposition}

\proof
We can assume that $X$ is connected. 
For any simplicial sheaf  $\mathcal F^\bullet$ on $X_\bullet$
consider the functor 
$a_*:\ker \mathcal F^0\stackrel{\delta_0^*-\delta_1^*}
{\longrightarrow}\mathcal F^1$
and let $Ra_*{\mathcal F}^\bullet$ be the total derived functor.
If $a_p:X_p\to X$ is the canonical map, we have the spectral sequence 
$$E_1^{p,q}=R^q(a_p)_* \mathcal F^p\Rightarrow R^{p+q}a_*{\mathcal F^\bullet}.$$
Since $H^i_\et(X_\bullet,\mathcal F^\bullet)= H^i_\et(X,Ra_*\mathcal F^\bullet)$,
it suffices to show that $a_*\Z\cong \Z$ and $R^1a_*\Z =0$.  
But $H^1_\et(Y,\Z)=0$ for normal schemes $Y$, so we have
$R^1(a_p)_* \Z=0$ for all $p$, and it suffices to show that
$$ 0\to \Z \to (a_0)_*\Z\to (a_1)_*\Z\to (a_2)_* \Z \to \cdots $$
is an exact sequence of sheaves on $X$. By the proper base-change
theorem for $H^0$ and constant sheaves \cite[II Remark 3.8]{milnebook}, 
we can assume that $X$
is the spectrum of a separably closed field. In this case, 
the hypercover has a section, and the sequence splits.
The second isomorphism
follows because both groups are calculated by the complex $\Z^{\pi_0(X_\bullet)}$
in degrees at most $1$. 

For the final statement, $H^1_\et(X,\Z)=0$ for normal $X$, hence
the exact sequence
$$ 0\to \Ext(H_{i-1}(D_{X_\bullet},\Z),\Z) \to H^i(D_{X_\bullet},\Z)\to 
\Hom(H_i(D_{X_\bullet},\Z),\Z) \to 0$$
and finite generation of $H_1(D_{X_\bullet},\Z)$ shows that the group is finite.
\proofend

\begin{lemma}\label{kklemma}
If $H_1(D_{X_\bullet},\Z)\otimes\Q_l/\Z_l=0$, then 
$H_1(\A_{X_\bullet}(k))\otimes \Q_l/\Z_l=0$.
\end{lemma}

\proof
Since the groups $H_i(D_{X_i},\Z)$ are finitely generated, 
the hypothesis implies 
that $H_1(D_{X_\bullet},\Z)$ is finite, and 
we obtain a short exact sequence 
$$ H_1(\A_{X_\bullet}^0(k),\Z)\to 
H_1(\A_{X_\bullet}(k),\Z)\to (finite) \to 0.$$
The result follows because $H_1(\A_{X_\bullet}^0(k),\Z)$ is an extension
of a finite group by the (divisible) group of $k$-rational
points of a semi-abelian variety, hence tensoring with $\Q_l/\Z_l$
annihilates it. 
\proofend

\noindent{\it Proof of Theorem \ref{rojtman}:}
Consider the following diagram, where the vertical maps
are induced by the albanese map:
\begin{equation}\begin{CD}\label{uioui}
H_1^S(X_\bullet,\Q_l/\Z_l)@>\sim>> H_0^S(X_\bullet,\Z)[l]
@>\sim>>H_0^S(X,\Z)[l]\\
@V\sim ValbV @VValbV @VValbV \\
H_1(\A_\bullet(k),\Q_l/\Z_l) @>\sim>> (\A_{X_0}(k)/d\A_{X_1}(k))[l] @>>>
\A_X(k)[l] 
\end{CD}\end{equation}
The left vertical map and the left horizontal maps are isomorphisms
by Theorem \ref{generalprop}, Lemma \ref{kklemma} and \eqref{ses3}, respectively.
Recall the exact sequence of presheaves
$$H_1(D_{X_\bullet},\Z)\stackrel{\delta}{\longrightarrow}
\A^0_{X_0} / d\A_{X_1}^0\to 
\A_{X_0}/d\A_{X_1}\to H_0(D_{X_\bullet},\Z) \to 0$$
from Lemma \ref{firstlemma}. The exactness of the global sections functor implies
that $\A_{X_0}(k)/d\A_{X_1}(k)\cong \big(\A_{X_0}/d\A_{X_1})(k)$,
and since $H_0(D_{X_\bullet},\Z)$ is torsion free, the right map in \eqref{uioui}
factors through $\big(\A_{X_0}^0/d\A_{X_1}^0+ \im \delta\big )(k)[l]$.
The finiteness of 
$H_1(D_{X_\bullet},Z)$ implies that $\A_{X_0}^0/d\A_{X_1}^0+ \im \delta$ 
is a semi-abelian variety, hence 
by Lemma \ref{firstlemma}, $\big(\A_{X_0}^0/d\A_{X_1}^0+ \im \delta\big )(k)[l]$
is the $l$-torsion of the connected
component of the largest locally semi-abelian scheme quotient of 
$\A_{X_0}/d\A_{X_1}$ which by Proposition \ref{albaba}
is the albanese scheme of $X$.
\proofend

\noindent{\it Proof of Theorem \ref{generalprop}, general case:} 
The only property of $X_\bullet$ used in the proof of
Theorem \ref{generalprop} are $H_1^S(X_i,\Z)\otimes \Q_l/\Z_l=0$ and 
$H_0^S(X_i,\Z)[l]\cong \A_{X_i}(k)[l]$.
By Theorem \ref{rojtman}, this hypothesis is satisfied 
for a hypercover by normal schemes. 
\proofend

\rem 
There is a canonical proper hypercover 
by normal schemes of a reduced semi-normal scheme $X$: 
Take $X_0=\tilde X$ to be the normalization of $X$. 
The diagonal map 
$\tilde X\to (\tilde X\times_X  \tilde X)^{\red} $ is a closed 
immersion with a 
section between reduced schemes, so is an irreducible component. 
The normalization of  $(\tilde X\times_X  \tilde X)^{\red}$ has 
$\tilde X$ as  a connected component, and we let $Z$ be its complement 
(of smaller dimension). Then
$$\tilde X\coprod Z  \to \tilde X \to X$$
is a $1$-truncated proper hypercover of $X$ by normal schemes.
Since the two projections from $\tilde X$ to itself are equal,
the difference of the maps induced on albanese schemes is trivial,
and the quotient of semi-abelian schemes is $\A_{\tilde X}(k)/d \A_{Z}(k)$
sitting in the exact sequence 
$$\ker \big( D_{Z}\to D_{\tilde X}\big) \to \A_{\tilde X}^0/d \A_{Z}^0(k)
\to \A_{\tilde X}(k)/d \A_{Z}(k)\to \Z^{\pi_0(X)}\to 0.$$ 
In particular, we obtain a surjection
$  {}_\tor H_0^S(X,\Z)\to  {}_\tor (\A_{\tilde X }(k)/d\A_{ Z}(k)) $.

\section{Curves}
We start with an example showing that the statement of the 
main theorem fails for non-normal curves.

Let $E$ be an elliptic curve and $p$ be
a closed point of $E$. Let $N$ be the variety obtained by glueing the points 
$0$ and $p$ of $E$. 
A proper hypercover of $N$ in low degree is given by 
$$ E\times_NE \times_NE\to  E\times_NE  
\stackrel{\delta_0,\delta_1}{\longrightarrow} E \to N.$$
The middle term is isomorphic to $E\cup x \cup y$ 
where $x$ and $y$ correspond to the points $(0,p)$ and $(p,0)$ in 
the product, respectively. Similarly, the term on the left is isomorphic
to $E$ and $6$ points corresponding to triples $(x,y,z)$ with $x,y,z\in \{0,p\}$
and not all equal. The albanese schemes are 
$$\begin{CD}
0@>>> E@>>> \A_2 @>>> \Z^7@>>> 0\\
@. @VVV @VVV @VVV \\
0@>>> E@>>> \A_1 @>>> \Z^3@>>>0\\
@.@V0VV @V\delta_1V\delta_0V @VVV\\
0@>>> E@>>> \A_0 @>>> \Z@>>> 0
\end{CD}$$ 
A calculation shows that $H_1(D_{X_\bullet},\Z)=\Z$, and the sequence
\eqref{ssse} becomes
$$H_1(\A_\bullet(k),\Z)\to  \Z\stackrel{\delta}{\to} E(k) \to 
H_0(\A_\bullet(k),\Z) \to \Z\to 0 ,$$
where $\delta$ sends $1$ to $p-0$ on $E$. 
Now assume that $p$ is not torsion. Then 
the albanese scheme of $N$ is the largest locally semi-abelian scheme quotient
of $\A_{X_0}$ modulo the subabelian variety generated by $\langle p \rangle $, 
hence it is isomorphic to $\Z$. 
The coranks of $H_1^S(N,\Q/\Z)$ and of ${}_\tor H_0^S(N,\Z)$ can be calculated to be
$3$, hence ${}_\tor H_0^S(N,\Z)$ is not isomorphic to the torsion of the 
albanese variety. 
However ${}_\tor H_0^S(N,\Z)$ isomorphic to the torsion of the abelian group 
quotient
$H_0(\A_\bullet(k),\Z) = 
\A_{X_0}(k)/d\A_{X_1}(k) \cong \big(\A_{X_0}^0/d\A_{X_1}^0\big) (k) /\im \delta$.
In other words, taking the quotient in the category of locally semi-abelian
schemes and then taking rational points does not give the correct answer,
but taking rational points, and then dividing in the category
of abelian groups, does.  More generally:

\begin{theorem}\label{curvecase}
Let $X$ be a reduced semi-normal curve. Then the albanese map
induces an isomorphism 
$$H_0^S(X,\Z) \cong  \A_{X_0}(k)/d\A_{X_1}(k)$$
\end{theorem}

By \eqref{ssse}, the right hand group is isomorphic to
$H_0(D_{X_\bullet},\Z)\oplus \A_{\tilde X}^0(k) /\im H_1(D_{X_\bullet},\Z)$, 
for $\tilde X$ the normalization of $X$, with $H_0(D_{X_\bullet},\Z)
\cong \Z^{\pi_0(X)}$ and 
$H_1(D_{X_\bullet},\Z)$ having the same rank as $H^1_\et(X,\Z)$.

\proof
We can choose $X_0$ to be the normalization of $X$, and $X_1=X_0 \coprod S$
with $S$ of dimension $0$. 
Since Suslin homology satisfies descent for hyperenvelopes  
\cite{ichdescent}, we obtain a commutative diagram with exact rows
\begin{equation}\label{curveproof}
\begin{CD}
H_0^S(X_1,\Z) @>>>   H_0^S(X_0,\Z) @>>>  H_0^S(X,\Z) @>>>  0\\
@V\alb_{X_1}VV @V\alb_{X_0}VV @V\alb_XVV  \\
\A_{X_1}(k) @>d>> \A_{X_0}(k)  @>>> \A_{X_0}(k)/d\A_{X_1}(k) @>>> 0.
\end{CD}
\end{equation}
The left two vertical maps are isomorphisms by the Abel-Jacobi theorem 
stating that for a regular curve the albanese map is an isomorphism. 
Hence the right hand map is an isomorphism. 
\proofend

{\bf Question:} Does the analog statement hold in higher dimensions,
i.e. is the surjection 
$$\alb_X: {}_\tor H_0^S(X,\Z) \to  {}_\tor (\A_{X_0}(k)/d\A_{X_1}(k))$$
an isomorphism for any reduced semi-normal scheme $X$?

The proof of Theorem \ref{curvecase} does not carry over, because there 
could be a uniquely divisible subgroup in
the albanese kernel of $X_0$ which maps to a torsion divisible group 
in the albanese kernel of $X$.  

On the other hand, the proof of Theorem \ref{curvecase} does carry
over if the albanese map is an isomorphism (not only on the torsion
part), as is known for the algebraic closure of finite fields,
and expected for the algebraic closure of the field of rational numbers.

\end{document}